\documentclass[12pt,leqno]{amsart}

\usepackage{soul}
\usepackage {yfonts}
\usepackage [T1] {fontenc}

\usepackage{xcolor}
\usepackage{textcomp}

\usepackage{enumitem}

\numberwithin{equation}{section}
\newtheorem{theorem}{Theorem}[section]
\newtheorem{prop}[theorem]{Proposition}
\newtheorem{lem}[theorem]{Lemma}

\newtheorem{rem}[theorem]{Remark}

\newcommand{\R}{\mathbb{R}}

\newcommand{\N}{\mathbb{N}}

\newcommand{\F}{\mathcal{F}}
\theoremstyle{plain}

\def\pn{\hfil\par\noindent}
\def\be{\begin{enumerate}}
\def\ds{\displaystyle }
\def\ee{\end{enumerate}}

\def\beqq{\begin{eqnarray*}}
\def\eeqq{\end{eqnarray*}}

\def\buildo#1\over#2{\mathrel{\mathop{\null#2}\limits^{#1}}}
\def\buildu#1\under#2{\mathrel{\mathop{\null#2}\limits_{#1}}}

\title{Asymptotic study of Leray Solution of 3D-NSE With Exponential Damping}
\author{Mongi Blel}
   \address{King Saud University, College of Sciences, Department of Mathematics,   Kingdom of Saudi Arabia}
   \email{mblel@ksu.edu.sa, jamelbenameur@gmail.com}
\author{Jamel Benameur}

\date{\today}


\begin{document}

\maketitle

\begin{abstract}
We study the uniqueness, the continuity in $L^2$ and the large
time decay for the Leray solutions of the $3D$ incompressible
Navier-Stokes equations with the nonlinear exponential
damping term $a (e^{b |u|^{\bf 2}}-1)u$, ($a,b>0$) studied by the second author in \cite{J1}.
\end{abstract}

\section{\bf Introduction}\ \\

 In this paper, we investigate the questions of the  existence, uniqueness and asymptotic study  of global weak solution to the following modified incompressible Navier-Stokes equations in $\R^3$

\begin{equation}\label{$S_2$}
 \left\{ \begin{matrix}
     \partial_t u
 -\nu\Delta u+ u.\nabla u  +a (e^{b |u|^2}-1)u = -\nabla p \hfill&\hbox{ in } \mathbb R^+\times \mathbb R^3\\
     {\rm div}\, u = 0 \hfill&\hbox{ in } \mathbb R^+\times \mathbb R^3\\
    u(0,x) =u^0(x) \;\;\hfill&\hbox{ in }\mathbb R^3,\\
    a,b>0\hfill&
\end{matrix}\right. \tag{$S$}
\end{equation}
 where $u=u(t,x)=(u_1,u_2,u_3)$, $p=p(t,x)$ denote respectively the unknown velocity and the unknown pressure of the fluid at the point $(t,x)\in \mathbb R^+\times \mathbb R^3$. The function $\nu$ denotes  the viscosity of fluid  and $u^0=(u_1^0(x),u_2^0(x),u_3^0(x))$ is the initial given velocity.
 The damping  of the system   comes from the resistance to the motion of the flow. It describes various physical situations such as porous media flow, drag or friction effects, and some dissipative mechanisms (see \cite{BD,BDC,H,HP} and references therein).
  The fact that ${\rm div}\,u = 0$, allows to write the term $(u.\nabla u):=u_1\partial_1 u+u_2\partial_2 u+u_3\partial_3u$ in the following form
$ {\rm div}\,(u\otimes u):=({\rm div}\,(u_1u),{\rm div}\,(u_2u),{\rm div}\,(u_3u)).$
 If the initial   velocity $u^0$ is quite regular, the divergence free condition determines uniquely the pressure $p$.\\
  Without loss of generality and in order to simplify the proofs of our results,  we consider the viscosity unitary ($\nu=1$).

  The global existence of weak solution of initial value problem of classical incompressible Navier-Stokes were proved by Leray and Hopf (see \cite{Hopf}-\cite{Leray}) long before. Uniqueness remains an open problem if the dimension $d\geq3$.\\
  The polynomial damping $\alpha|u|^{\beta-1}u$ is studied in \cite{CJ} by Cai and Jiu. They proved the global  existence of weak solution in
 $$L^\infty(\R^+,L^2(\R^3))\cap L^2(\R^+,\dot H^1(\R^3))\cap L^{\beta+1}(\R^+,L^{\beta+1}(\R^3)).$$
 The exponential damping $a (e^{b |u|^2}-1)u$ is studied in \cite{J1}  by J. Benameur. He   proved the existence of global weak solution in
 $$L^\infty(\R^+,L^2(\R^3))\cap L^2(\R^+,\dot H^1(\R^3))\cap \mathcal E_b,$$
 where
 $\ds \mathcal E_b=\{f:\R^+\times\R^3\rightarrow\R\ :{\rm measurable}, \ (e^{b|f|^2}-1)|f|^2\in L^1(\R^+\times\R^3)\}.$\\

 The purpose of this paper is to prove the uniqueness and the continuity of the global   solution given in \cite{J1}. Using   Friedritch method, we construct  approximate solutions and we make more delicate estimates to proceed the compactness arguments. In particular, we obtain some new a priory estimates:
$$\|u(t)\|_{L^2}^2+2\int_0^t\|\nabla u(s)\|_{L^2}^2ds+2 a\int_0^t\|(e^{b |u(s)|^2}-1)|u(s)|^2\|_{L^1}ds\leq \|u^0\|_{L^2}^2,$$

  To prove the uniqueness we use the energy method and the approximates systems. The proof of the asymptotic study is based on a decomposition of the solution in high and low frequencies and the uniqueness of such  solution in a well chosen time $t_0$.

In our case of exponential damping, we find more regularity of Leray
solution in $\cap_pL^p(\R^+,L^p(\R^3))$. In particular, we give a new energy estimate. Our main result is the following:

\begin{theorem}\label{th1}\pn
 Let $u^0\in L^2(\mathbb R^3)$ be a divergence free vector fields, then there is a unique global solution of the system  $(S)$:
$u\in C_b(\R^+,L^2(\mathbb R^3))\cap L^2(\R^+,\dot H^1(\mathbb
R^3))\cap\mathcal E_b$. Moreover, for all $t\geq0$
\begin{equation}\label{eqth1-1}
\|u(t)\|_{L^2}^2+2\int_0^t\|\nabla
u(s)\|_{L^2}^2ds+2a\int_0^t\|(e^{b
|u(s)|^2}-1)|u(s)|^2\|_{L^1}ds\leq \|u^0\|_{L^2}^2.
\end{equation}
Moreover, we have
\begin{equation}\label{eqth1-2}
\limsup_{t\to \infty}\|u(t)\|_{L^2}=0.
\end{equation}
\end{theorem}

\begin{rem}\pn
\begin{enumerate}
\item The new results in this theorem is the uniqueness, the
continuity of the global week   solution in   $L^2(\R^3)$    and its asymptotic behavior at infinity.

\item The uniqueness of weak solution implies that
\begin{equation}\label{eqth1-1}
\|u(t_2)\|_{L^2}^2+2\int_{t_1}^{t_2}\|\nabla
u(s)\|_{L^2}^2ds+2a\int_{t_1}^{t_2}\|(e^{b
|u(s)|^2}-1)|u(s)|^2\|_{L^1}ds\leq \|u(t_1)\|_{L^2}^2,
\end{equation}
which implies that $(t\rightarrow\|u(t)\|_{L^2})$ is decreasing.
\end{enumerate}
\end{rem}

\section{\bf Notations and Preliminary Results}
The
Friedrich operator $J_R$  is defined for $R>0$ by:
  $\ds J_R(D)f=\F^{-1}(\chi_{B_R} \widehat{f}),$
where $B_R$ is the ball of center $0$ and radius $R$ and  $f\in L^2(\R^3)$.
The  Leray operator  $\mathbb P$ is the projector operator of $ (L^2(\R^3))^3$ on the   space of divergence-free vector fields $L^2_\sigma(\R^3)$.\\
If $f$ is in the  Schwartz space $ (\mathcal S(\R^n))^3$,
$$\mathcal F(\mathbb P f)=\widehat{f}(\xi)-(\widehat{f}(\xi).\frac{\xi}{|\xi|})\frac{\xi}{|\xi|}=M(\xi)\widehat{f}(\xi) $$
and
 $\ds \left(\mathbb P  f\right)_k(x) = \frac{1}{(2\pi)^{\frac 3 2}} \int_{\R^3}  \left( \delta_{kj}-\frac{\xi_k \xi_j}{ \vert \xi \vert^2}\right)
\widehat{  f}_j(\xi) \, e^{i \xi \cdot  x}\,   d\xi,$ where
  $M(\xi)$ is  the matrix $(\delta_{k,\ell}-\frac{\xi_k\xi_\ell}{|\xi|^2})_{1\leq k,\ell\leq 3} $. \\

Define also the operator $A_R(D)$ on $L^2(\R^3)$ by:
 $$\ds A_R(D)u=\mathbb P J_R(D)u=\mathcal F^{-1}(M(\xi)\chi_{B_R}(\xi)\widehat{u}).$$
%


 To simplify the exposition of the main result, we first collect some preliminary results  and we give some new technical lemmas.

\begin{prop}(\cite{HBAF})\label{prop1}\pn
 Let $H$ be a Hilbert space.
\begin{enumerate}
\item The unit ball is weakly compact, that is: if $(x_n)$ is a bounded sequence   in $H$, then there is a subsequence $(x_{\varphi(n)})$ such that
$$(x_{\varphi(n)}|y)\to  (x|y),\;\forall y\in H.$$

\item If $x\in H$ and $(x_n)$   a bounded sequence   in $H$ such that
$\ds\lim_{n\to+\infty}(x_n|y)=  (x|y)$, for all $y\in H,$
then $\|x\|\leq\ds \liminf_{n\to \infty}\|x_n\|.$

\item If $x\in H$ and $(x_n)$ is a bounded sequence   in $H$ such that\\
$\ds\lim_{n\to+\infty}(x_n|y)=  (x|y)$, for all $y\in H$
and
 $\limsup_{n\to \infty}\|x_n\|\leq \|x\|,$
then $\ds \lim_{n\to \infty}\|x_n-x\|=0.$
\end{enumerate}
\end{prop}
We recall   the following product law in the homogeneous Sobolev spaces:

\begin{lem}(\cite{JYC})\label{lem1}\pn
Let $s_1,\ s_2$ be two real numbers and $d\in\N$.

\begin{enumerate}
\item If $s_1<\frac d 2$\; and\; $s_1+s_2>0$, there exists a constant  $C_1=C_1(d,s_1,s_2)$, such that: if $f,g\in \dot{H}^{s_1}(\mathbb{R}^d)\cap \dot{H}^{s_2}(\mathbb{R}^d)$, then $f.g \in \dot{H}^{s_1+s_2-\frac{d}{2}}(\mathbb{R}^d)$ and
$$\|fg\|_{\dot{H}^{s_1+s_2-\frac{d}{2}}}\leq C_1 (\|f\|_{\dot{H}^{s_1}}\|g\|_{\dot{H}^{s_2}}+\|f\|_{\dot{H}^{s_2}}\|g\|_{\dot{H}^{s_1}}).$$

\item If $s_1,s_2<\frac d 2$\; and\; $s_1+s_2>0$ there exists a constant $C_2=C_2(d,s_1,s_2)$ such that: if $f \in \dot{H}^{s_1}(\mathbb{R}^d)$\; and\; $g\in\dot{H}^{s_2}(\mathbb{R}^d)$, then  $f.g \in \dot{H}^{s_1+s_2-\frac{d}{2}}(\mathbb{R}^d)$ and
$$\|fg\|_{\dot{H}^{s_1+s_2-\frac{d}{2}}}\leq C_2 \|f\|_{\dot{H}^{s_1}}\|g\|_{\dot{H}^{s_2}}.$$
\end{enumerate}
 \end{lem}

\begin{lem}\label{lem2}(\cite{JL})\pn
    Let $ \beta>0$ and $d\in\N$. Then, for all $x,y\in\R^d$, we have

    \begin{equation}\label{eqn-lem2-1}
    \langle |x|^{\beta}x-|y|^{\beta}y ,x-y\rangle\geq \frac{1}{2}(|x|^{\beta}+|y|^{\beta})|x-y|^{2},
    \end{equation}
and
 \small
\begin{equation}\label{eqn-lem2-2}
\langle (e^{b |x|^2}-1)x-(e^{b |y|^2}-1)y,x-y\rangle\geq \frac{1}{2}\Big((e^{b |x|^2}-1)+(e^{b |y|^2}-1)\Big)|x-y|^{2}.
\end{equation}
\end{lem}
\begin{prop}\label{prop2}(\cite{J1})\pn
Let $\nu_1,\nu_2,\nu_3\in[0,\infty)$, $r_1,r_2,r_3\in(0,\infty)$ and $f^0\in L^2_\sigma(\R^3)$. \\
For $n\in\N$, let $F_n:\R^+\times\R^3\to \R^3$ be a measurable function in $C^1(\R^+,L^2(\R^3))$ such that $$A_n(D)F_n=F_n,\;F_n(0,x) =A_n(D)f^0(x)$$ and

\begin{enumerate}
\item [(E1)]
$\ds  \partial_t F_n+\sum_{k=1}^3\nu_k|D_k|^{2r_k} F_n+ A_n(D){\rm div}\,(F_n\otimes F_n)+ A_n(D)h(|F_n|)F_n =0.$

\item [(E2)]
\beqq
&&\ds \|F_n(t,.)\|_{L^2}^2+2\sum_{k=1}^3\nu_k\int_0^t\||D_k|^{r_k} F_n(s,.)\|_{L^2}^2ds\\
&&\hskip 2cm +2 a\int_0^t\|h(|F_n(s,.)|)|F_n(s,.)|^2\|_{L^1}ds \leq \|f^0\|_{L^2}^2.
\eeqq
\end{enumerate}
 where  $\ds h(\lambda)= a(e^{b \lambda^2}-1),$\, with $a,b  >0$.
  Then: for every $\varepsilon>0$ there is $\delta=\delta(\varepsilon,a,b,\nu_1,\nu_2,\nu_3,r_1,r_2,r_3,\|f^0\|_{L^2})>0$
  such that: for all $t_1,t_2\in\R^+$, we have

\begin{equation}\label{eqn-1}
\Big(|t_2-t_1|<\delta\Longrightarrow \|F_n(t_2)-F_n(t_1)\|_{H^{-s_0}}<\varepsilon\Big),\;\forall n\in\N,
\end{equation}
with $\ds  s_0\ge \max(3,2r_1,2r_2,2r_3).$
\end{prop}
\begin{lem}\label{lemf1} Let $a,b>0$, then there is a unique real $\lambda_0=\lambda_0(a,b)\geq0$ such that: For all $\lambda\geq0$
$$a(e^{b\lambda}-1)\leq \lambda\Longrightarrow \lambda\in[0,\lambda_0].$$
Precisely\\
$\bullet$ If $ab\geq1$, we have $\lambda_0=0$,\\
$\bullet$ If $ab<1$, we have $\lambda_0>\frac{1}{b}\log(1/ab)$.
\end{lem}

\section{\bf Proof of the Main Theorem \ref{th1}}\ \\
The proof of the theorem is given in four steps:

\subsection{Existence of  Week Solution}\ \\
In this step, we build approximate solutions of the system $(S)$ inspired by the method used in  \cite{J1,JYC}, hence we construct a global solution. For this, consider the approximate system with   parameter $n\in\N$:
$$(S_n)
  \begin{cases}
     \partial_t u
 -\Delta J_nu+ J_n(J_nu.\nabla J_nu)  + a J_n[(e^{b |J_nu|^2}-1)J_nu] =\;\;-\nabla p_n\hbox{ in } \mathbb R^+\times \mathbb R^3\\
 p_n=(-\Delta)^{-1}\Big({\rm div}\,J_n(J_nu.\nabla J_nu)  + a {\rm div}\,J_n[(e^{b |J_nu|^2}-1)J_nu]\Big)\\
     {\rm div}\, u = 0 \hbox{ in } \mathbb R^+\times \mathbb R^3\\
    u(0,x) =J_nu^0(x) \;\;\hbox{ in }\mathbb R^3.
  \end{cases}
$$
  $J_n$ is the Friedritch operator defined in the second section.

\begin{enumerate}
\item[$\bullet$] By Cauchy-Lipschitz Theorem, we obtain a unique solution $u_n\in C^1(\R^+,L^2_\sigma(\R^3))$ of $(S_{2,n})$. Moreover, $J_nu_n=u_n$ such that

\begin{equation}\label{eqn-6}
\|u_n(t)\|_{L^2}^2+2\int_0^t\|\nabla
u_n\|_{L^2}^2+2a\int_0^t\|(e^{b
|u_n|^2}-1)|u_n|^2\|_{L^1}\leq \|u^0\|_{L^2}^2.
\end{equation}

\item[$\bullet$]
The sequence $(u_n)_n$ is bounded in $L^\infty(\R^+,L^2(\R^3))$ and on $L^2(\R^+,\dot H^{1}(\R^3)$. Using Proposition \ref{prop2} and the  interpolation method, we deduce that the sequence $(u_n)_n$ is equicontinuous on $H^{-1}(\R^3)$.

\item[$\bullet$] Let $(T_q)_q $ be a strictly  increasing sequence such that  $\ds\lim_{q\to+\infty} T_q=\infty$. Consider a sequence
of functions $(\theta_q)_{q }$ in $C_0^\infty(\R^3)$  such that

$$\left\{\begin{array}{l}
\theta_q(x)=1,\ {\rm for}\  |x|\le   q+\frac{5}{4}\\
\theta_q(x)=0,\ {\rm for}\   |x|\ge  q+2 \\
0\leq \theta_q\leq 1.
\end{array}\right.$$
Using  \eqref{eqn-6}, the equicontinuity of the sequence $(u_n)_n$   on $H^{-1}(\R^3)$
 and classical argument by combining Ascoli's theorem and the Cantor diagonal process, there exists a subsequence $(u_{\varphi(n)})_n$   and\\
$u\in L^\infty(\R^+,L^2(\R^3))\cap C(\R^+,H^{-3}(\R^3))$ such that: for all $q\in\N$,

\begin{equation}\label{eqn-7}
\lim_{n\to \infty}\|\theta_q(u_{\varphi(n)}-u)\|_{L^\infty([0,T_q],H^{-4})}=0.
\end{equation}
In particular, the sequence $(u_{\varphi(n)}(t))_n$ converges weakly in $L^2(\R^3)$ to $u(t)$ for all
$t\geq0$.

\item[$\bullet$] Combining the above inequalities, we obtain:

\begin{equation}\label{eqn-8}
\|u(t)\|_{L^2}^2\!+\!2\int_0^t\!\|\nabla
u(s)\|_{L^2}^2ds\!+\!2a\int_0^t\!\|(e^{b |u(s)|^2} -1)
|u(s)|^2\|_{L^1}ds\!\leq\! \|u^0\|_{L^2}^2.
\end{equation}
for all $t\geq0$.

\item[$\bullet$] $u$ is a solution of the system $(S)$.\end{enumerate}

\subsection{Continuity of the Solution in $L^2$}\pn
  In this section, we give a simple proof of the continuity of the solution $u$ of the system $(S)$ and we prove also that $u\in C(\R^+,L^2(\R^3))$. The construction of the solution is based on the Friedrich approximation method. We point out that we can use this method to show the same results as in  \cite{HP}.\\
$\bullet$ By inequality (\ref{eqn-8}), we get
$$\limsup_{t\to 0}\|u(t)\|_{L^2}\leq\|u^0\|_{L^2}.$$
Then, proposition \ref{prop1}-(3) implies that

$$\limsup_{t\to 0}\|u(t)-u^0\|_{L^2}=0,$$
which ensures the continuity of $u$ at $0$.\\
$\bullet$ Consider the  functions
$$v_{n,\varepsilon}(t,.)=u_{\varphi(n)}(t+\varepsilon,.),\;p_{n,\varepsilon}(t,.)=p_{\varphi(n)}(t+\varepsilon,.),$$
for $n\in\N$ and $\varepsilon>0$. We have:

$$\begin{array}{lcl}
\partial_tu_{\varphi(n)}-\Delta u_{\varphi(n)}+J_{\varphi(n)}(u_{\varphi(n)}.\nabla
u_{\varphi(n)})+a J_{\varphi(n)}(e^{b|u_{\varphi(n)}|^2}-1)u_{\varphi(n)}&=&-\nabla p_{\varphi(n)} \\
\partial_tv_{n,\varepsilon}-\Delta v_{n,\varepsilon}+J_{\varphi(n)}
(v_{n,\varepsilon}.\nabla v_{n,\varepsilon})+a J_{\varphi(n)}(e^{b|v_{n,\varepsilon}|^4}-1)v_{n,\varepsilon}&=&-\nabla p_{n,\varepsilon}\\
\end{array}$$
The function $w_{n,\varepsilon}=u_{\varphi(n)}-v_{n,\varepsilon}$ fulfills the following:

\beqq
&&\partial_tw_{n,\varepsilon}-\Delta w_{n,\varepsilon}  +a
J_{\varphi(n)}\Big((e^{b|u_{\varphi(n)}|^2}-1)u_{\varphi(n)}-(e^{b|v_{n,\varepsilon}|^2}-1)v_{n,\varepsilon}\Big)\\
&&\hskip  3cm = -\nabla (p_{\varphi(n)}-p_{n,\varepsilon})+J_{\varphi(n)}(w_{n,\varepsilon}.\nabla w_{n,\varepsilon})\\
&&\hskip  3cm-J_{\varphi(n)}(w_{n,\varepsilon}.\nabla u_{\varphi(n)})
 - J_{\varphi(n)}(u_{\varphi(n)}.\nabla w_{n,\varepsilon}).
 \eeqq
Taking the scalar product in $L^2(\R^3)$ with
$w_{n,\varepsilon}$ and using the properties ${\rm div}\ w_{n,\varepsilon}=0$ and
$\langle w_{n,\varepsilon}.\nabla w_{n,\varepsilon},w_{n,\varepsilon}\rangle=0$, we get

\begin{eqnarray}\label{eqn-9}
\frac{1}{2}\frac{d}{dt}\|w_{n,\varepsilon}\|_{L^2}^2+\|\nabla
w_{n,\varepsilon}\|_{L^2}^2& +&a \langle
J_{\varphi(n)}\Big((e^{b|u_{\varphi(n)}|^2}-1)u_{\varphi(n)}
-(e^{b|v_{n,\varepsilon}|^4}-1)v_{n,\varepsilon}\Big);w_{n,\varepsilon}\rangle_{L^2}
\nonumber\\
&=& -\langle J_{\varphi(n)}(w_{n,\varepsilon}.\nabla u_{\varphi(n)});w_{n,\varepsilon}\rangle _{L^2}.
\end{eqnarray}
Using inequality  \eqref{eqn-lem2-2}, we get

\beqq &&\langle
J_{\varphi(n)}\Big((e^{b|u_{\varphi(n)}|^2}-1)u_{\varphi(n)}-(e^{b|v_{n,\varepsilon}|^2}-1)v_{n,\varepsilon}\Big);
w_{n,\varepsilon}\rangle _{L^2}\\
 &&\hskip 5cm=\langle (e^{b|u_{\varphi(n)}|^2}-1)u_{\varphi(n)}-
 (e^{b|v_{n,\varepsilon}|^2}-1)v_{n,\varepsilon};J_{\varphi(n)}w_{n,\varepsilon}\rangle _{L^2}\\
 &&\hskip 5cm= \langle (e^{b|u_{\varphi(n)}|^2}-1)u_{\varphi(n)}-(e^{b|v_{n,\varepsilon}|^2}-1)v_{n,\varepsilon};w_{n,
\varepsilon}\rangle_{L^2}\\
 && \hskip 5cm\geq  \frac{1}{2}\int_{\R^3}\Big((e^{b|u_{\varphi(n)}|^2}-1)+(e^{b|v_{n,\varepsilon}|^2}-1)\Big)|w_{n,\varepsilon}|^2\\
&&\hskip 5cm  \geq   \frac{1}{2}\int_{\R^3}(e^{b|u_{\varphi(n)}|^2}-1)|w_{n,\varepsilon}|^2,
\eeqq

$$\begin{array}{lcl}
|\langle J_{\varphi(n)}(w_{n,\varepsilon}.\nabla u_{\varphi(n)});
w_{n,\varepsilon}\rangle _{L^2}|&\leq&\ds \int_{\R^3}|w_{n,\varepsilon}|.|u_{\varphi(n)}|.|\nabla w_{n,\varepsilon}|\\
&\leq&\ds \frac{1}{2}\int_{\R^3}|w_{n,\varepsilon}|^2|u_{\varphi(n)}|^2+\frac{1}{2}\|\nabla w_{n,\varepsilon}\|_{L^2}^2.
\end{array}$$

Combining the  identity  \eqref{eqn-lem2-2} and the inequality \eqref{eqn-9}, we get

$$\frac{d}{dt}\|w_{n,\varepsilon}\|_{L^2}^2+\|\nabla w_{n,\varepsilon}\|_{L^2}^2+a\int_{\R^3}(e^{b|u_{\varphi(n)}|^2}-1)|w_{n,\varepsilon}|^2\leq
\int_{\R^3}|w_{n,\varepsilon}|^2|u_{\varphi(n)}|^2.$$

By using Lemma \ref{lemf1}, we get:\\
If $ab\geq1$
$$\frac{d}{dt}\|w_{n,\varepsilon}\|_{L^2}^2+\|\nabla w_{n,\varepsilon}\|_{L^2}^2\leq0,$$
If $ab<1$, we have
$$\frac{d}{dt}\|w_{n,\varepsilon}\|_{L^2}^2+\|\nabla w_{n,\varepsilon}\|_{L^2}^2+a\int_{\R^3}(e^{b|u_{\varphi(n)}|^2}-1)|w_{n,\varepsilon}|^2\leq
+\int_{A_t}|w_{n,\varepsilon}|^2|u_{\varphi(n)}|^2,$$
where $A_{n,t}=\{x\in\R^3/\;a(e^{b|u_{\varphi(n)}(t,x)|^2}-1)<|u_{\varphi(n)}(t,x)|^2\}$. Then, also by Lemma \ref{lemf1} we obtain
$$x\in A_{n,t}\Longrightarrow |u_{\varphi(n)}(t,x)|^2\leq \lambda_0,$$
which implies
$$\frac{d}{dt}\|w_{n,\varepsilon}\|_{L^2}^2+\|\nabla w_{n,\varepsilon}\|_{L^2}^2\leq
\lambda_0\int_{A_t}|w_{n,\varepsilon}|^2\leq \lambda_0\||w_{n,\varepsilon}\|_{L^2}^2.$$
In all cases, we get
$$\frac{d}{dt}\|w_{n,\varepsilon}\|_{L^2}^2\leq \lambda_0\||w_{n,\varepsilon}\|_{L^2}^2.$$
By Gronwall Lemma, we get
$$\|w_{n,\varepsilon}(t)\|_{L^2}^2\leq  \|w_{n,\varepsilon}(0)\|_{L^2}^2e^{\lambda_0t}.$$
Then
$$\|u_{\varphi(n)}(t+\varepsilon)-u_{\varphi(n)}(t)\|_{L^2}^2\leq
\|u_{\varphi(n)}(\varepsilon)-u_{\varphi(n)}(0)\|_{L^2}^2e^{\lambda_0t
}.$$ For $t_0>0$ and $\varepsilon\in(0,t_0)$, we have

$$\|u_{\varphi(n)}(t_0+\varepsilon)-u_{\varphi(n)}(t_0)\|^{2}_{L^{2}}
\leq\|u_{\varphi(n)}(\varepsilon)-u_{\varphi(n)}(0)\|^{2}_{L^{2}}
\exp(2\lambda_0t_0),$$

$$\|u_{\varphi(n)}(t_0-\varepsilon)-u_{\varphi(n)}(t_0)\|^{2}_{L^{2}}\leq
\|u_{\varphi(n)}(\varepsilon)-u_{\varphi(n)}(0)\|^{2}_{L^{2}}
\exp(2\lambda_0t_0).$$
So
\beqq
\| u_{\varphi(n)}(\varepsilon)-u_{\varphi(n)}(0) \|_{L^2}^2 &=&
\| J_{\varphi(n)} u_{\varphi(n)}(\varepsilon)-J_{\varphi(n)}u_{\varphi(n)}(0) \|_{L^2}^2 \\
&=&\|\chi_{\varphi(n)}\left(\widehat{u_{\varphi(n)}}  -\widehat{u^0}\right)\|_{\varphi(n)}^2\\
&\le&\| u_{\varphi(n)}(\varepsilon)-u^0 \|_{L^2}^2\\
&\le& 2\|u^0 \|_{L^2}^2-2Re\langle
u_{\varphi(n)}(\varepsilon),u^0\rangle. \eeqq But $\ds
\lim_{n\to+\infty}\langle
u_{\varphi(n)}(\varepsilon),u^0\rangle=\langle u
(\varepsilon),u^0\rangle$. Hence
$$\liminf_{n\to \infty}\|u_{\varphi(n)}(\varepsilon)-u_{\varphi(n)}(0)\|^{2}_{L^{2}}\leq 2\|u^0\|^{2}_{L^{2}}-
2Re\langle u(\varepsilon);u^0\rangle_{L^2}.$$ Moreover, for all
$q,N\in\N$ \beqq
\|J_N\Big(\theta_q.(u_{\varphi(n)}(t_0\pm\varepsilon)-u_{\varphi(n)}(t_0))\Big)\|^{2}_{L^2}
 &\leq & \|\theta_q.(u_{\varphi(n)}(t_0\pm\varepsilon)-u_{\varphi(n)}(t_0))\|^{2}_{L^2}\\
 &\leq&
 \|u_{\varphi(n)}(t_0\pm\varepsilon)-u_{\varphi(n)}(t_0)\|^{2}_{L^2}.
 \eeqq

 Using  \eqref{eqn-7} we get, for $q$ big enough,
 $$\|J_N\Big(\theta_q.(u(t_0\pm\varepsilon)-u(t_0))\Big)\|^{2}_{L^2}
 \leq \liminf_{n\to \infty}\|u_{\varphi(n)}(t_0\pm\varepsilon)-u_{\varphi(n)}(t_0)\|^{2}_{L^2}.$$
 Then
$$\|J_N\Big(\theta_q.(u(t_0\pm\varepsilon)-u(t_0))\Big)\|^{2}_{L^2}
 \leq 2\Big(\|u^0\|^{2}_{L^{2}}-Re\langle u(\varepsilon);u^0\rangle_{L^2}\Big)\exp(2\lambda_0t_0).$$
By applying the monotone convergence theorem in the order $N\to \infty$ and $q\to \infty$, we get
$$\|u(t_0\pm\varepsilon,.)-u(t_0,.))\|^{2}_{L^2}
 \leq 2\Big(\|u^0\|^{2}_{L^{2}}-Re\langle u(\varepsilon);u^0\rangle_{L^2}\Big)\exp(2\lambda_0t_0).$$
Using the continuity at 0 and make $\varepsilon\to 0$, we get the continuity at $t_0$.
\subsection{Uniqueness of the Solution}\ \\
Let $u,v$ be two solutions of $(S)$ in the space

$$C_b(\R^+,L^2(\R^3))\cap L^2(\R^+,\dot H^1(\R^3))\cap \mathcal E_b.$$
The function $w=u-v$ satisfies the following:

$$\partial_tw-\Delta w+a \Big((e^{b|u|^2}-1)u-(e^{b|v|^2}-1)v\Big)= -\nabla (p-\tilde p)+w.\nabla w-w.\nabla u-
u.\nabla w  .$$
Taking the scalar product in $L^2$ with $w$, we get

$$\frac{1}{2}\frac{d}{dt}\|w\|_{L^2}^2+\|\nabla w\|_{L^2}^2+a \langle \Big((e^{b|u|^2}-1)u-(e^{b|v|^2}-1)v\Big);w\rangle _{L^2}=-\langle w.\nabla u;w\rangle _{L^2}  .$$
The idea is to lower the term $\langle \Big((e^{b|u|^2}-1)u-(e^{b|v|^2}-1)v\Big);w\rangle _{L^2}$ with the help of the Lemma \ref{lem2} and then divide the term find into two equal pieces, one to absorb the nonlinear term and the other is used in the last inequality.\\

By using inequality \eqref{eqn-lem2-2}, we get

\beqq
\langle \Big((e^{b|u|^2}-1)u-(e^{b|v|^2}-1)v\Big);w\rangle _{L^2}&\geq & \frac{1}{2}\int_{\R^3}  \Big((e^{b|u|^4}-1)+(e^{b|v|^4}-1)\Big)|w|^2\\
&\geq& \frac{1}{2}\int_{\R^3} (e^{b|u|^4}-1)|w|^2.
\eeqq
Moreover, we have
\beqq
|\langle w.\nabla u;w\rangle _{L^2}|&=&|\langle {\rm div}\,(w\otimes u);w\rangle _{L^2}|=|\langle w\otimes u;\nabla w\rangle _{L^2}|\\
&\leq&\ds \int_{\R^3}|w|.|u|.|\nabla w|
\leq\ds \frac{1}{2}\int_{\R^3}|w|^2|u|^2+\frac{1}{2}\|\nabla w\|_{L^2}^2\\
\eeqq
Combining the above inequalities, we get
$$\frac{d}{dt}\|w\|_{L^2}^2+\|\nabla w\|_{L^2}^2+a\int_{\R^3}(e^{b|u|^2}-1)|w|^2\leq \int_{\R^3}|w|^2|u|^2.$$
By using Lemma \ref{lemf1} and the set
$$A_t=\{x\in\R^3/\;a(e^{b|u|^2}-1)<|u(t,x)|^2\},$$
we get
$$\frac{d}{dt}\|w\|_{L^2}^2\leq \lambda_0\|w\|_{L^2}^2.$$
and, Gronwall Lemma gives
$$\|w\|_{L^2}^2\leq \|w^0\|_{L^2}^2e^{\lambda_0t}.$$
As $w^0=0$, then $w=0$ and $u=v$. Which implies the uniqueness.
\subsection{Asymptotic Study of the Global Solution}\pn
In this subsection we  prove the asymptotic behavior (\ref{eqth1-2}). Here, we use a modified version of Benameur-Selmi method(\cite{JR}). The idea is to write the nonlinear term of exponential type as follows
$$a(e^{b|u|^2}-1)u=a(e^{b|u|^2}-1-b|u|^2)u+a|u|^2u.$$
The first term, one treats it of the matter of \cite{BB1}, and the second and by specifying the difficulty in small frequencies, we descend into the Sobolev norm of negative index $H^{-\sigma}(\mathbb R^3)$,\,$\sigma>0$.\\
For this, let $\varepsilon>0$ be positif real number. For $\delta>0$, put the following functions
$$v_\delta=\mathcal F^{-1}({\bf 1}_{B(0,\delta)}(\xi)\widehat{u}(\xi)),\;w_\delta=u-v_\delta.$$
We have
$$v_\delta=\sum_{k=1}^4f_{\delta,k}(t),$$
where
$$\begin{array}{lcl}
f_{\delta,1}&=&e^{t\Delta}v_\delta^0,\;\;\;v_\delta^0=\mathcal F^{-1}({\bf 1}_{B(0,\delta)}(\xi)\widehat{u^0}(\xi))\\
f_{\delta,2}&=&-\displaystyle\int_0^te^{(t-z)\Delta}{\bf 1}_{B(0,\delta)}(D)\mathbb P{\rm div}\,(u\otimes u)\\
f_{\delta,3}&=&-\displaystyle a\int_0^te^{(t-z)\Delta}{\bf 1}_{B(0,\delta)}(D)\mathbb P(e^{b|u|^2}-1-b|u|^2)u\\
f_{\delta,4}&=&-\displaystyle\int_0^te^{(t-z)\Delta}{\bf 1}_{B(0,\delta)}(D)\mathbb P(|u|^2u).
\end{array}$$
$\bullet$ We have
$$\|f_{\delta,1}(t)\|_{L^2}\leq\|v_\delta^0\|_{L^2}.$$
As $\lim_{\delta\rightarrow0}\|v_\delta^0\|_{L^2}=0$, then there is $\delta_1>0$ such that
\begin{equation}\label{eqjm1}
\sup_{t\geq0}\|f_{\delta,1}(t)\|_{L^2}<\varepsilon/8,\;\forall 0<\delta<\delta_1.
\end{equation}
$\bullet$ We have
$$\begin{array}{lcl}
\|f_{\delta,2}(t)\|_{H^{-1/4}}&\leq&\displaystyle\int_0^t\|e^{(t-z)\Delta}{\bf 1}_{B(0,\delta)}(D)\mathbb P{\rm div}\,(u\otimes u)\|_{H^{-1/4}}\\
&\leq&\displaystyle\int_0^t\|{\bf 1}_{B(0,\delta)}{\rm div}\,(u\otimes u)\|_{H^{-1/4}}\\
&\leq&\displaystyle\int_0^t\|{\bf 1}_{B(0,\delta)}|D|\,(u\otimes u)\|_{H^{-1/4}}\\
&\leq&\displaystyle\int_0^t\|{\bf 1}_{B(0,\delta)}|D|\,(u\otimes u)\|_{H^{-1/4}}\\
&\leq&\displaystyle\delta^{1/4}\int_0^t\|{\bf 1}_{B(0,\delta)}|D|^{3/4}\,(u\otimes u)\|_{\dot H^{-1/4}}\\
&\leq&\displaystyle\delta^{1/4}\int_0^t\||D|^{3/4}\,(u\otimes u)\|_{\dot H^{-1/4}}\\
&\leq&\displaystyle\delta^{1/4}\int_0^t\|(u\otimes u)\|_{\dot H^{1/2}}\\
&\leq&\displaystyle C\delta^{1/4}\int_0^t\|u\|_{\dot H^{1}}^2,\;(s_1+s_2=2,s_1=s_2=1)\\
&\leq&\displaystyle C\|u^0\|_{L^2}^2\delta^{1/4}.
\end{array}$$
But $$\|f_{\delta,2}(t)\|_{L^2}\leq c_0(1+\delta^2)^{1/4}\|f_{\delta,2}(t)\|_{H^{-1/4}}$$
Then
$$\|f_{\delta,2}(t)\|_{L^2}\leq c_0C\|u^0\|_{L^2}^2(1+\delta^2)^{1/4}\delta^{1/4}.$$
there is $\delta_2>0$ such that
\begin{equation}\label{eqjm2}
\sup_{t\geq0}\|f_{\delta,2}(t)\|_{L^2}<\varepsilon/8,\;\forall 0<\delta<\delta_2.
\end{equation}
$\bullet$ We have
$$\begin{array}{lcl}
\|f_{\delta,3}(t)\|_{H^{-2}}&\leq&\displaystyle\int_0^t\|e^{(t-z)\Delta}{\bf 1}_{B(0,\delta)}(D)\mathbb P(e^{b|u|^2}-1-b|u|^2)u\|_{H^{-2}}\\
&\leq&\displaystyle\int_0^t\|{\bf 1}_{B(0,\delta)}(D)(e^{b|u|^2}-1-b|u|^2)u\|_{H^{-2}}\\
&\leq&\displaystyle R_\delta\int_0^t\|(e^{b|u|^2}-1-b|u|^2)u\|_{L^1},
\end{array}$$
where
$$R_\delta=\|{\bf 1}_{B(0,\delta)}(D)\|_{H^{-2}}=\Big(\int_{B(0,\delta)}\frac{1}{(1+|\xi|^2)^2}d\xi\Big)^{1/2}
\leq\Big(\int_{B(0,\delta)}d\xi\Big)^{1/2}=c_0\delta^{3/2}.$$
By using the elementary inequality with $M_{b}>0$
$$(e^{bz^2}-1-bz^2)z\leq M_{b}(e^{bz^2}-1)z^2,\;\forall z\geq0$$
we get
$$\int_0^t\|(e^{b|u|^2}-1-b|u|^2)u\|_{L^1}\leq M_b\int_0^t\|(e^{b|u|^2}-1)|u|^2\|_{L^1}\leq (2a)^{-1}M_b\|u^0\|_{L^2}^2.$$
Combining the above inequalities we get
$$\|f_{\delta,3}(t)\|_{H^{-2}}\leq c_0(2a)^{-1}M_b\|u^0\|_{L^2}^2\delta^{3/2}.$$
But $$\|f_{\delta,3}(t)\|_{L^2}\leq c_0(1+\delta^2)^{2}\|f_{\delta,2}(t)\|_{H^{-2}}$$
then
$$\|f_{\delta,3}(t)\|_{L^2}\leq c_0(2a)^{-1}M_b\|u^0\|_{L^2}^2(1+\delta^2)^{2}\delta^{3/2}.$$

Then there is $\delta_3>0$ such that
\begin{equation}\label{eqjm3}
a\sup_{t\geq0}\|f_{\delta,3}(t)\|_{L^2}<\varepsilon/8,\;\forall 0<\delta<\delta_3.
\end{equation}
$\bullet$ By using Lemma \ref{lem1} with the well choice of $s_1$ and $s_2$, we get
$$\begin{array}{lcl}
\|f_{\delta,4}(t)\|_{H^{-1/2}}&\leq&\displaystyle\int_0^t\|e^{(t-z)\Delta}{\bf 1}_{B(0,\delta)}(D)\mathbb P|u|^2u\|_{H^{-1/2}}\\
&\leq&\displaystyle\int_0^t\|{\bf 1}_{B(0,\delta)}(D)|u|^2u\|_{H^{-1/2}}\\
&\leq&\displaystyle\int_0^t\|{\bf 1}_{B(0,\delta)}(D)|u|^2u\|_{\dot H^{-1/2}}\\
&\leq&\displaystyle\int_0^t\|{\bf 1}_{B(0,\delta)}(D)|D|^{1/2}|D|^{-1/2}|u|^2u\|_{\dot H^{-1/2}}\\
&\leq&\displaystyle\delta^{1/2}\int_0^t\|{\bf 1}_{B(0,\delta)}(D)|D|^{-1/2}|u|^2u\|_{\dot H^{-1/2}}\\
&\leq&\displaystyle\delta^{1/2}\int_0^t\||D|^{-1/2}|u|^2u\|_{\dot H^{-1/2}}\\
&\leq&\displaystyle\delta^{1/2}\int_0^t\||u|^2u\|_{\dot H^{-1}}\\
&\leq&\displaystyle\delta^{1/2}C\int_0^t\|u\|_{\dot H^{1/2}}\||u|^2\|_{L^2},\;(s_1+s_2=1/2,s_1=0,\;s_2=1/2)\\
&\leq&\displaystyle\delta^{1/2}C'\int_0^t\|u\|_{\dot H^{1/2}}^2\|u\|_{\dot H^1},\;(s_1+s_2=3/2,\;s_1=1/2,\;s_2=1),\\
&\leq&\displaystyle\delta^{1/2}C'\int_0^t\|u\|_{L^2}\|u\|_{\dot H^1}^2,\;(by\;interpolation),\\
&\leq&\displaystyle\delta^{1/2}C''\|u^0\|_{L^2}\int_0^t\|\nabla u\|_{L^2}^2\\
&\leq&\displaystyle\delta^{1/2}C''\|u^0\|_{L^2}^3.
\end{array}$$
But $$\|f_{\delta,4}(t)\|_{L^2}\leq c_0(1+\delta^2)^{1/2}\|f_{\delta,2}(t)\|_{H^{-1/2}}$$
then
$$\|f_{\delta,4}(t)\|_{L^2}\leq c_0C''\|u^0\|_{L^2}^3(1+\delta^2)^{1/2}\delta^{1/2}.$$
Then there is $\delta_4>0$ such that
\begin{equation}\label{eqjm4}
\sup_{t\geq0}\|f_{\delta,4}(t)\|_{L^2}<\varepsilon/8,\;\forall 0<\delta<\delta_4.
\end{equation}
Combining the above equations (\ref{eqjm1})-(\ref{eqjm2})-(\ref{eqjm3})-(\ref{eqjm4}) we get
\begin{equation}\label{eqjm5}
\sup_{t\geq0}\|v_{\delta_0}(t)\|_{L^2}<\varepsilon/2,\;\delta_0=\frac{1}{2}\min_{1\leq i\leq 4}\delta_i.
\end{equation}

In other hand we have
$$\begin{array}{lcl}
\displaystyle\int_0^\infty\|w_{\delta_0}(t)\|_{L^2}^2dt&\leq&\displaystyle\delta_0^{-2}\int_0^\infty\|\nabla w_{\delta_0}(t)\|_{L^2}^2dt\\
&\leq&\displaystyle\delta_0^{-2}\int_0^\infty\|\nabla w_{\delta_0}(t)\|_{L^2}^2dt\\
&\leq&\displaystyle\delta_0^{-2}\|u^0\|_{L^2}^2\\
\end{array}$$
As $(t\rightarrow\|w_{\delta_0}(t)\|_{L^2}^2)$ is continuous, then there is a time $t_0\geq0$ such that
\begin{equation}\label{eqjm6}
\|w_{\delta_0}(t_0)\|_{L^2}^2<\varepsilon/2.
\end{equation}
Combining inequalities (\ref{eqjm5})-(\ref{eqjm6}) we get
$$\|u(t_0)\|_{L^2}\leq \|v_{\delta_0}(t_0)\|_{L^2}+\|w_{\delta_0}(t_0)\|_{L^2}<\varepsilon.$$
As $(t\rightarrow\|w_{\delta_0}(t)\|_{L^2}^2)$ is decreasing, then
$$\|u(t)\|_{L^2}<\varepsilon,\;\forall t\geq t_0.$$
which complete the proof.

\end{document}